\documentclass[10pt]{article}
\usepackage[latin1]{inputenc}
\usepackage{amsmath}
\usepackage{amsfonts}
\usepackage{amssymb}
\usepackage[american]{babel}
\usepackage{graphicx}
\usepackage{amsthm}
\usepackage{enumerate}

\begin{document}

\title{Covering a cubic graph with perfect matchings}

\author{G.Mazzuoccolo \footnote{Dipartimento di Scienze Fisiche, Informatiche e Matematiche, Universit\`a di Modena e Reggio Emilia, via Campi 213/b, 41125 Modena (Italy)}}

\maketitle

\newtheorem{proposition}{Proposition}[section]
\newtheorem{theorem}{Theorem}[section]
\newtheorem{lemma}{Lemma}[section]
\newtheorem{definition}{Definition}[section]
\newtheorem{example}{Example}[section]
\newtheorem{corollary}{Corollary}[section]
\newtheorem{conjecture}{Conjecture}[section]

\begin{abstract} \noindent 
Let $G$ be a bridgeless cubic graph. A well-known conjecture of Berge and Fulkerson can be stated as follows: there exist five perfect matchings of $G$ such that each edge of $G$ is contained in at least one of them. Here, we prove that in each bridgeless cubic graph there exist five perfect matchings covering a portion of the edges at least equal to $\frac{215}{231}$. By a generalization of this result, we decrease the best known upper bound, expressed in terms of the size of the graph, for the number of perfect matchings needed to cover the edge-set of $G$. 
\end{abstract}

\noindent \textit{ Keywords: Berge-Fulkerson conjecture, perfect matchings, cubic graphs. }

\section{Introduction}\label{sec:intro}

Throughout this paper, a graph $G$ always means a simple connected finite graph (without loops and parallel edges). Furthermore, along the entire paper $G$ stands for a bridgeless cubic graph, unless otherwise specified, and we denote by $V(G)$ and $E(G)$ the vertex-set and the edge-set of $G$, respectively. A perfect matching of $G$ is a $1$-regular spanning subgraph of $G$. Following the definition introduced in \cite{KaiKraNor} and \cite{Pat}, we define $m_t(G)$ to be the maximum fraction of the edges in $G$ that can be covered by $t$ perfect matchings, and by $m_t$ the infimum of all $m_t(G)$ over all bridgeless cubic graphs. That is,

$$m_t=\inf_{G} \max_{M_1,\ldots,M_t} \frac{|\bigcup_{i=1}^t M_i|}{|E(G)|}$$

The Berge--Fulkerson conjecture, one of the challenging open problems in graph theory, can be easily stated in terms of $m_t$ : it is the assertion that $m_5=1$ (see \cite{Maz}).

Kaiser, Kr\'al and Norine \cite{KaiKraNor} proved that $m_2=\frac{3}{5}$ and $m_3 \geq \frac{27}{35}$. In a final remark of their paper, they announced, without a proof, the more general inequality $a_t \leq m_t$, where $a_t$ is the sequence defined by the recurrence 
\begin{equation}\label{a_t}
a_t=\frac{t}{2t+1}(1-a_{t-1})+a_{t-1} 
\end{equation}
and $a_0=0$. 

In the present paper we present a complete proof of the result announced in \cite{KaiKraNor}, and we use it to deduce a new upper bound in terms of $t$ for the size of a bridgeless cubic graph admitting a covering with $t$ perfect matchings. More precisely, it was known (see for instance \cite{Pat}) that a bridgeless cubic graph with fewer than $(\frac{3}{2})^t$ edges can be covered using $t$ perfect matchings. We improve this bound by proving that each bridgeless cubic graph with fewer than $\frac{2^t}{\sqrt t}$ edges can be covered with $t$ perfect matchings.

Finally, as a by-product of our main result, we also obtain that in each bridgeless cubic graph there is a set of $t$ perfect matchings with no $(2t+1)$-cut in their intersection, giving partial support to a conjecture of Kaiser and Raspaud \cite{KaiRas}, in a special form due to M\'acajov\'a
and Skoviera \cite{MacSko}, about the existence of two perfect matchings with no odd cut in their intersection.

The main tool for our proof is the Perfect Matching Polytope Theorem of Edmonds (see \cite{Edm}); we briefly recall it in the next section.

\section{The perfect matching polytope}
 
Let $G$ be a graph. A \textit{minimal cut} $C$ in $G$ is a subset of $E(G)$ such that $G \setminus C$ has more components than $G$ does, and $C$ is inclusion-wise minimal with this property. A $k$-cut is a minimal cut of cardinality $k$. When $X \subseteq V(G)$, let $\partial X$ denote the set of edges with precisely one end in $X$.

Let $w$ be a vector in $\mathbb{R}^{E(G)}$. The entry of $w$ corresponding  to an edge  $e$ is denoted by $w(e)$, and for $A \subseteq E(G)$, we define the \textit{weight} $w(A)$ of $A$ as $\sum_{e \in A} w(e)$. The vector $w$ is a \textit{fractional perfect matching} of $G$ if it satisfies the following properties:
\begin{enumerate}[a)]
 \item $0 \leq w(e) \leq 1$ for each $e \in E(G)$,
 \item $w(\partial\{v\})=1$ for each vertex $v \in V$, and
 \item $w(\partial X) \geq 1$ for each $X \subseteq V(G)$ of odd cardinality.
\end{enumerate}

We will denote by $P(G)$ the set of all fractional perfect matchings of $G$. If $M$ is a perfect matching, then the characteristic vector $\chi^M \in \mathbb{R}^E$ of $M$ is contained in $P(G)$. Furthermore, if $w_1,\ldots,w_n \in P(G)$, then any convex combination $\sum_{i=1}^n \alpha_i w_i$ also belongs to $P(G)$. It follows that $P(G)$ contains the convex hull of all vectors $\chi^M$ such that $M$ is a perfect matching of $G$. The Perfect Matching Polytope Theorem of Edmonds asserts that the converse inclusion also holds:

\begin{theorem}[Edmonds]
For any graph $G$, the set $P(G)$ is precisely the convex hull of the characteristic vectors of perfect matchings of $G$.
\end{theorem}

The main tool in our proof is the following property of a fractional perfect matching:

\begin{lemma}\label{tool}
If $w$ is a fractional perfect matching in a graph $G$ and $c \in \mathbb{R}^E$, then $G$ has a perfect matching $M$ such that
$$ c \cdot \chi^M \geq c \cdot w $$
where $\cdot$ denotes the dot product. Moreover, there exists such a perfect matching $M$ that contains exactly one edge of each odd cut $C$ with $w(C)=1$.  
\end{lemma}

A proof of this lemma can be found in \cite{KaiKraNor}.
 
\section{A lower bound for $m_t$}\label{sec:lowerbound}

In this section we prove $m_t \geq a_t$ for each index $t$, where $a_t$ is the sequence defined by the recurrence 
\begin{equation}\label{a_t}
a_t=\frac{t}{2t+1}(1-a_{t-1})+a_{t-1} 
\end{equation}
and $a_0=0$. In particular, we deduce $m_5 \geq \frac{215}{231}$.\\

Let ${\cal M}^t=\{M_1,\ldots,M_{t}\}$, where each $M_i$ is a perfect matching of a bridgeless cubic graph $G$, and set ${\cal M}^0=\emptyset$. 

For each subset $A$ of the edge-set of $G$ we define $$\Phi(A,{\cal M}^t)=\sum_{i=1}^{t} |A \cap M_i|. $$ 

Furthermore, define the weight $w_{{\cal M}^t}$ by

$$ w_{{\cal M}^t}(e)=\frac{t+1-\sum_{i=1}^{t} |M_i \cap \{e\}|}{2t+3}=\frac{t+1- \Phi(\{e\},{\cal M}^t)}{2t+3}.$$

When $|{\cal M}|=1$ or $2$, we obtain the spaecial cases $w_2$ and $w_3$ in the main theorem of \cite{KaiKraNor}. 

Hence, for a set $A \subseteq E(G)$, with $|A|=k$, the weight of $A$ is given by the following relation: 

$$ w_{{\cal M}^t}(A)=\frac{k(t+1)-\Phi(A,{\cal M}^t)}{2t+3}.$$

In other words, by the definition of the weight $w_{{\cal M}^t}$, the weight of a set depends only on the size of the intersections with the perfect matchings in ${\cal M}^t$.

An easy calculation proves the following lemma,

\begin{lemma}\label{Phi}
For $A \subseteq E(G)$ with $|A|=k$. 
$$ w_{{\cal M}^t}(A) \geq 1  \iff \Phi(A,{\cal M}^t) \leq t(k-2)+(k-3)$$
and equality holds on one side of the implication if and only if it holds on the other side.
\end{lemma}

We are now ready to state the main result.

\begin{theorem}\label{maintheorem}
Let $a_t$ be the sequence (\ref{a_t}).
Then, $m_{t} \geq a_{t}$ for each index $t$.  
\end{theorem}
\textit{Proof.} We argue by induction: let ${\cal M}^t=\{M_1,\ldots,M_t\}$ be a set of $t$ perfect matchings such that $w_{{\cal M}^t}$ is a fractional perfect matching of a bridgeless cubic graph $G$. Now, we construct a set ${\cal M}^{t+1}$ of $t+1$ perfect matchings and such that  $w_{{\cal M}^{t+1}}$ is a fractional perfect matching of $G$. Set $c_t=1-\chi^{\bigcup_{i=1}^{t}M_i}$. By Lemma \ref{tool} there exists a perfect matching $M_{t+1}$ such that 
$$ c_t \cdot \chi^{M_{t+1}} \geq c_t \cdot w_{{\cal M}^t} $$ and such that $M_{t+1}$ contains exactly one edge for each cut $C$ with  $w_{{\cal M}^t}(C)=1$.

In order to prove that $w_{{\cal M}^{t+1}}$, where ${\cal M}^{t+1}={\cal M}^{t} \cup \{ M_{t+1} \}$, is a fractional perfect matching of $G$, we have to verify the properties $(a),(b)$ and $(c)$ of the definition of fractional perfect matching.\\
$(a)$ $0\leq w_{{\cal M}^{t+1}}(e) \leq 1$ trivially holds for each $e \in E(G)$. \\
$(b)$ Let $v$ be a vertex of $G$. We have 
$$ w_{{\cal M}^{t+1}}(\partial \{v\})=\frac{3(t+2)-\Phi(\partial \{v\},{\cal M}^{t+1})}{2t+5}$$
since each perfect matching $M_i$ intersects $\partial \{v\}$ exactly one time we have $\Phi(\partial \{v\},{\cal M}^{t+1})=t+1$ and then 
$$ w_{{\cal M}^{t+1}}(\partial \{v\})=\frac{3(t+2)- (t+1)}{2t+5}=1.$$
$(c)$ Let $C$ be a $k$-cut of $G$, with $k$ odd.

$$ w_{{\cal M}^{t+1}}(C)=\frac{k(t+2)-\Phi(C,{\cal M}^{t+1})}{2t+5}$$

If $k=3$, then by the assumption that $w_{{\cal M}^t}$ is a fractional perfect matching,  we have 
$$w_{{\cal M}^{t}}(C)=\frac{3(t+1)-\Phi(C,{\cal M}^{t})}{2t+3} \geq 1;$$
that is, $\Phi(C,{\cal M}^{t}) \leq t = |{\cal M}^t|$.
Furthermore, a $3$-cut intersects each of the $t$ perfect matchings $M_i$, hence $\Phi(C,{\cal M}^{t}) \geq t$. Then $\Phi(C,{\cal M}^{t})=t$ and $w_{{\cal M}^{t}}(C)=1$. By Lemma \ref{tool}, we obtain  $|C \cap M_{t+1}|=1$. Now, we can compute $\Phi(C,{\cal M}^{t+1})=   
\Phi(C,{\cal M}^{t})+|C \cap M_{t+1}|=t+1$. We have proved that
$$ w_{{\cal M}^{t+1}}(C)=\frac{3(t+2)-\Phi(C,{\cal M}^{t+1})}{2t+5}=\frac{3(t+2)-(t+1)}{2t+5}=1,$$
as required.

If $k>3$, we distinguish two cases according that  $w_{{\cal M}^{t}}(C)=1$ (that is $\Phi(C,{\cal M}^{t})=t(k-2)+(k-3)$) or  $w_{{\cal M}^{t}}(C)>1$ (that is $\Phi(C,{\cal M}^{t})<t(k-2)+(k-3)$). 

In the former case we have $|C \cap M_{t+1}|=1$, so
$$\Phi(C,{\cal M}^{t+1})=\Phi(C,{\cal M}^{t})+|C \cap M_{t+1}|=t(k-2)+(k-3)+1=(t+1)(k-2)\leq(t+1)(k-2)+(k-3).$$
Now  $w_{{\cal M}^{t+1}}(C) \geq 1$ by Lemma \ref{Phi}.\\ 
 
In the latter case we have $|M_{t+1} \cap C|\leq k$, hence 

$$\Phi(C,{\cal M}^{t+1})=\Phi(C,{\cal M}^{t})+|C \cap M_{t+1}| < t(k-2)+(k-3)+k.$$

Since each perfect matching intersects an odd cut an odd number of times  follows $\Phi(C,{\cal M}^{t}) \equiv t \pmod 2$, so the previous inequality is equivalent to
$$\Phi(C,{\cal M}^{t+1}) \leq t(k-2)+(k-5)+k = (t+1)(k-2)+(k-3).$$
By Lemma \ref{Phi}, we can infer $w_{{\cal M}^{t+1}}(C) \geq 1$ also in this case.

By induction, since the basic step $w_{{\cal M}^0}(e)=\frac{1}{3}$ is trivially a fractional perfect matching, we have that $w_{{\cal M}^t}$ (with ${\cal M}^t$ constructed as described above) is a fractional perfect matching for each value of $t$. Therefore, the following holds:

$$ c_t \cdot \chi^{M_{t+1}} \geq c_t \cdot w_{{\cal M}^t}.$$
The left side of the previous inequality is exactly the number of edges of $M_{t+1}$ not covered by ${\cal M}^t$, while the right side is $\frac{t}{2t+1}$ times the number of edges not covered by ${\cal M}^t$. Denoting by $a_t$ the fraction of edges of $E(G)$ covered by ${\cal M}^t$, we obtain
$$ m_{t+1} - a_{t} \geq a_{t+1}-a_{t} = (1-a_t) \cdot \frac{t+1}{2t+3},$$
$$ m_{t+1} \geq (1-a_t) \cdot \frac{t+1}{2t+3}+a_t = a_{t+1}$$
and the assertion follows. \qed \\

By  direct calculation, $a_5=\frac{215}{231}$, and the following corollary holds: 

\begin{corollary}
Let $G$ be a bridgeless cubic graph. There exist five perfect matchings of $G$ that cover at least $\lceil \frac{215}{231}|E(G)| \rceil$ edges of $G$. 
\end{corollary}

Furthermore, we obtain the following corollary by the proof of Theorem \ref{maintheorem}.

\begin{corollary}
Let $G$ be a bridgeless cubic graph and $t$ a positive integer. There exist $t$ perfect matchings of $G$ with no $(2t+1)$-cut in their intersection. 
\end{corollary}

\section{A bound for the size of a bridgeless cubic graph which admits a covering with $t$ perfect matchings}
As remarked in \cite{Pat}, it is unknown whether $m_t=1$ for any $t\geq5$. The best known result in this direction is the following: if $G$ is a bridgeless cubic graph of size $|E(G)|$, then $m_t=1$ when $t>log_{\frac{3}{2}}(E)$. In other words the edge-set of each bridgeless cubic graph with fewer than $(\frac{3}{2})^t$ edges can be covered by $t$ perfect matchings. 
The following improvement of that bound is a consequence of Theorem \ref{maintheorem}.

\begin{theorem}
If $G$ is a bridgeless cubic graph with fewer than $\frac{2^t}{\sqrt{t}}$ edges, then there is a covering of $G$ by $t$ perfect matchings. 
\end{theorem}
\textit{Proof.} Let $G$ be a bridgeless cubic graph of size $|E(G)|$. It is trivial that if $|E(G)| m_t > |E(G)|-1$, that is $|E(G)| < \frac{1}{1-m_t}$, then there exists a covering of $G$ by $t$ perfect matchings.
The inequality $ \frac{1}{1-a_t} \leq \frac{1}{1-m_t}$ follows by Theorem \ref{maintheorem}. 
Therefore the theorem is proved, if we prove that the inequality $\lfloor \frac{2^t}{\sqrt{t}} \rfloor < \frac{1}{1-a_t}$ holds for each $t>0$.
One can easily verify that the inequality holds for $t<5$.
Now we will argue by induction to prove that $\frac{2^t}{\sqrt{t}} < \frac{1}{1-a_t}$ holds for each $t\geq5$.
For $t=5$, we have by direct computation  $$\frac{32}{\sqrt{5}}<\frac{1}{1-a_5}=\frac{231}{16}.$$ 
Suppose $\frac{2^t}{\sqrt{t}}< \frac{1}{1-a_{t}}$. We have
$$a_{t+1}=\frac{t+1}{2t+3}(1-a_t)+a_t>\frac{t+1}{2t+3}+\frac{t+2}{2t+3}(1-\frac{\sqrt{t}}{2^t})=1-\frac{t+2}{2t+3}\frac{\sqrt{t}}{2^t};$$
the assertion follows by a direct check that $$\frac{t+2}{2t+3}\frac{\sqrt{t}}{2^t}<\frac{\sqrt{t+1}}{2^{t+1}}$$ holds for each $t$.\qed
        
\section{Final Remarks}

We would like to stress that the special cases $m_2 \geq a_2 = \frac{3}{5}$ and $m_3 \geq a_3= \frac{27}{35}$ of Theorem \ref{maintheorem} were already proved in \cite{KaiKraNor}. Furthermore, it can be immediately checked that the Petersen graph realizes the equality $m_2=a_2$. Mainly for this reason, it is conjectured in \cite{Pat} that also $m_3$ and $m_4$ reach their minimum when $G$ is the Petersen graph.

It is trivial that every possible counterexample to the Berge--Fulkerson conjecture must be a snark, if it exists: a stronger version of the Berge--Fulkerson conjecture could be that $m_4(G)=1$ for each snark other than the Petersen graph (see \cite{FouVan}). Essentially the Petersen graph should be the unique obstruction to obtaining a covering with at most four perfect matchings\footnote{After this paper was submitted, a first counterexample for this conjecture was constructed in \cite{Hag} and  an infinite family of examples was presented in \cite{EspMaz}.}. Following this point of view, Bonvicini and the author prove in \cite{BonMaz} that a large class of bridgeless cubic graphs (including the Petersen graph) of order $2n$ can be covered by at most $4$ matchings of size $n-1$.\\
Another possible future study could be to prove an analogous of Theorem \ref{maintheorem} for a special class of $r$-regular graph defined in \cite{Sey}: Seymour proposed a natural generalization of the Berge-Fulkerson conjecture for this class of graphs and it turns out that his conjecture is equivalent to the assertation that each graph $G$ of this class can be covered with at most $2r-1$ perfect matchings (see \cite{Maz2}). It remains completely open the problem of establishing which portion of the edges of a $r$-regular graph, for $r>3$, can be covered with $t$ perfect matchings. 

\section{Acknowledgments}
The author is grateful to the reviewers and the editor for their very useful comments and advices.

\end{document}